\documentstyle[12pt]{amsart}

\catcode`\@=11

\long\def\@savemarbox#1#2{\global\setbox#1\vtop{\hsize\marginparwidth 
  \@parboxrestore\tiny\raggedright #2}}
\newcommand\lref[1]{\ref{#1}%
\@ifundefined{r@DisplaY #1}{}{ (#1)}}

\newcommand\fakelabel[2]{\@bsphack\if@filesw {\let\thepage\relax
   \newcommand\protect{\noexpand\noexpand\noexpand}%
\xdef\@gtempa{\write\@auxout{\string
      \newlabel{#1}{{#2}{\thepage}}}}}\@gtempa
   \if@nobreak \ifvmode\nobreak\fi\fi\fi\@esphack}

\catcode`\@=12

\def\Empty{}
\newcommand\oplabel[1]{
  \def\OpArg{#1} \ifx \OpArg\Empty {} \else
        \label{#1}
  \fi}
\newtheorem{theoremSt}{Theorem}[section]

\newtheorem{exampleSt}[theoremSt]{Example}
\newtheorem{exerciseSt}[theoremSt]{Exercise}

\newcommand\MakeStEnv[1]{
  \newenvironment{#1}[1]{
  \begin{#1St} \oplabel{##1}%
  \global\def\CrntSt{\thetheoremSt}%
}{ 
  \end{#1St} }
  \newenvironment{#1+}[1]{
  \begin{#1St} \label{##1}%
  \label{DisplaY ##1}%
  \global\def\CrntSt{\thetheoremSt}%
  \def\Labl{##1}\ifx\Labl\Empty{} \else {\em (\Labl)\,}\fi%
}{ 
  \end{#1St} }
}
\MakeStEnv{theorem}
\MakeStEnv{corollary}
\MakeStEnv{proposition}
\MakeStEnv{lemma}
\MakeStEnv{definition}
\MakeStEnv{conjecture}

\newlength{\saveu}

\newcommand{\startproof}[1]{%
\medbreak\mbox{}\noindent{\it Proof of #1:}%
}
\newcommand{\finishproof}[1]{ 
  \def\FPArg{#1}
  \ifx\FPArg\Empty
        \newcommand\FPArg{\CrntSt}  \fi
  \smallbreak\noindent\makebox[\textwidth]{\hfill\fbox{\FPArg}}
  \medbreak\noindent
}

\newcommand\FF{{\cal F}}

\newcommand\LL{{\cal L}}
\newcommand\MM{{\cal M}}

\newcommand\PP{{\cal P}}

\newcommand\PMF{{\PP\kern-2pt\MM\FF}}

\newcommand\PML{{\PP\kern-2pt\MM\LL}}

\newcommand\half{{\textstyle{1\over2}}}

\newcommand\ep{\epsilon}

\newcommand\union{\cup}
\newcommand\intersect{\cap}
\newcommand\bbR{{\mathord{\text{I\kern-2pt R}}}}        
\newcommand\bbH{{\mathord{\text{I\kern-2pt H}}}}        

\newcommand\C{{\bold C}}

\newcommand\Z{{\bold Z}}
\newcommand\R{{\bold R}}

\newcommand\Hyp{{\bold H}}

\newcommand\PSL[1]{\text{PSL}_{#1}}

\newcommand\bigrightarrow[1]{\hbox to #1{\rightarrowfill}}
\newcommand\bigleftarrow[1]{\hbox to #1{\leftarrowfill}}

\newcommand\boundary{\partial}
\newcommand\semidir{\mathrel{\hbox{\vrule depth-.03ex height1.1ex\kern-0.15em$\times$}}}

\newcommand\til{\widetilde}

\newcommand{\diam}{\operatorname{diam}}

\numberwithin{equation}{section}

\begin{document}

\title[Spectral theory, Hausdorff dimension and topology]{Spectral
theory, Hausdorff dimension and the topology of hyperbolic 3-manifolds}  
\author{Richard D. Canary}
\address{Department of Mathematics, University of Michigan,
Ann Arbor, MI 48109}
\author{Yair N. Minsky}
\address{Department of Mathematics, SUNY-Stony Brook, Stony Brook, NY 11794}
\author{Edward C. Taylor}
\address{Department of Mathematics, University of Michigan,
Ann Arbor, MI 48109}

\date{August 12, 1996}

\thanks{To appear in J. Geometric Analysis.
First author partially supported by a grant from the National
Science Foundation and a fellowship from the Alfred P. Sloan
Foundation.  Second author partially supported by an NSF Postdoctoral
Fellowship.  Third author partially supported by a Rackham
fellowship.}

\maketitle

\begin{abstract}
Let $M$ be a compact 3-manifold whose interior admits a complete hyperbolic
structure. We let $\Lambda(M)$ be the supremum of $\lambda_0(N)$ where
$N$ varies over all hyperbolic 3-manifolds homeomorphic to the
interior of $M$. Similarly, we let $D(M)$ be the infimum of the Hausdorff
dimensions of limit sets of  Kleinian groups whose quotients are
homeomorphic to the interior of $M$. We observe that $\Lambda(M)=D(M)(2-D(M))$
if $M$ is not handlebody or a thickened torus.
We characterize exactly when $\Lambda(M)=1$ and $D(M)=1$
in terms of the characteristic submanifold of the incompressible core
of $M$.
\end{abstract}

\section{Introduction}

When a closed 3-manifold admits a hyperbolic structure, this structure
is unique by Mostow's rigidity theorem \cite{mostow}. It follows that any
invariant of the hyperbolic structure is automatically a topological
invariant. 
One example of this is the hyperbolic volume, which agrees with 
Gromov's simplicial norm (see \cite{gromov}).

In this paper we will consider geometrically derived invariants for
compact 3-manifolds with boundary whose interiors admit complete
hyperbolic metrics of infinite volume. In this case the hyperbolic
structure is not unique, and in fact, Thurston's geometrization
theorem together with the Ahlfors-Bers quasiconformal deformation
theory guarantee that there is at least a one complex dimensional
space of hyperbolic structures.

To obtain a topological invariant in this context, one may 
begin with a natural geometric invariant of a hyperbolic
metric, and minimize (or maximize) it over the
class of all hyperbolic metrics on a given 3-manifold.
In particular, given a hyperbolic 3-manifold
$N=\Hyp^3/\Gamma$ we will consider the bottom $\lambda_0(N)$ of the
$L^2$ spectrum of the Laplacian, and the Hausdorff dimension $d(N)$ of
the limit set $L_\Gamma$ of $\Gamma$. (See section 2 for more precise
definitions).  Although these two geometric
invariants seem very different, the work of Patterson, Sullivan and
others has established that they are closely related.

Call a compact, orientable 3-manifold $M$
{\em hyperbolizable} if its interior
admits a complete hyperbolic metric. This is a topological condition:
Thurston's geometrization theorem 
asserts that an orientable compact 3-manifold with non-empty
boundary is hyperbolizable
if and only if it is irreducible and atoroidal. Given a hyperbolizable
$M$, 
we will define 
$$\Lambda(M) = \sup \lambda_0(N)$$
and 
$$D(M) = \inf d(N),$$
where $N$ varies over the space of all hyperbolic 3-manifolds homeomorphic
to the interior of $M$.

Very roughly speaking, we will find that $D$ increases with topological
complexity, and $\Lambda$ decreases. 
The purpose of this paper is to
establish some quantitative aspects of this intuition.
In particular we will characterize exactly which manifolds have
$\Lambda(M)=1$, the maximal possible value. 
 
Let us introduce some topological notation. A compact irreducible
manifold $M$ has an {\em incompressible core}, which is a (possibly
disconnected) submanifold
with incompressible boundary, from which $M$ is obtained by
adding 1-handles. A compact irreducible
3-manifold $M$ with incompressible boundary is called a {\em generalized
book of $I$-bundles} if one may find a disjoint collection $A$ of
essential annuli in $M$ such that each component $R$ of the manifold
obtained by cutting $M$ along $A$ is either a solid torus, a thickened
torus, or homeomorphic to an $I$-bundle such that $\partial R\cap
\partial M$ is the associated $\partial I$-bundle.

The following is a ``scorecard'' for the basic facts about $\Lambda$
and $D$ for hyperbolizable manifolds.

\medskip

\begin{center}
\small
\begin{tabular}{|p{1.5in}|c|c|l|} \hline
$M$ (hyperbolizable) & $\Lambda$ & $D$ & Relation \\ \hline
Handlebody or thickened torus & 1 & 0 & \\ \hline
Not solid or thickened torus, but $\boundary M$ is a union of tori  & 0 & 2 & 
        \parbox[t]{1in}{\vspace{5ex} $\Lambda = D(2-D)$} \\ \cline{1-3}
Incompressible core consists of
generalized books of $I$-bundles & 1 & 1 & \\ \cline{1-3}
Any other manifold with boundary & $0<\Lambda<1$ & $1<D<2$ & \\ \hline
\end{tabular}
\end{center}

%

\medskip\par\noindent
Our main theorems fill in the first two entries in the last two rows. 
The rest of the entries follow from known results, as explained in
section 2.

\medskip
\par\noindent
{\bf Main Theorem I:} {\em Let $M$ be a compact, orientable,
hyperbolizable 3-manifold.
$\Lambda(M)=1$ if and only if every component of the incompressible
core of $M$ is a generalized book of $I$-bundles. Otherwise,
$\Lambda(M) <1$.}

\medskip
\par\noindent
{\bf Main Theorem II:} {\em Let $M$ be a compact, orientable, hyperbolizable 
3-manifold which is not a handlebody or a thickened torus.
$D(M)=1$ if and only if every component of its incompressible
core is a generalized book of $I$-bundles. Otherwise, $D(M)>1$.}

\medskip

One could think of these results as analogues of the fact that the Gromov
norm of a closed, irreducible 3-manifold  $M$ is zero
if and only if there exists
a collection $T$ of incompressible tori in $M$ such that each component of
$M-T$ is a Seifert fibre space. In particular, if $M$ has incompressible
boundary, then $D(M)=1$ if and only if the Gromov norm of the double of
$M$ is 0.

\section{Preliminaries}
In this section we will more carefully define our invariants, derive
their basic properties and summarize the proof of the main theorems.

\subsection{Definitions}
If $N$ is a complete orientable hyperbolic 3-manifold, then it is isometric
to the quotient of  hyperbolic 3-space $\Hyp^3$ by a group $\Gamma$ of
orientation-preserving isometries. Any orientation-preserving
isometry extends continuously to  a conformal transformation of the
sphere at infinity $S^2_\infty$ of hyperbolic 3-space. We define
the {\em domain of discontinuity}
$\Omega_\Gamma$ to be the maximal open subset of $S^2_\infty$
on which $\Gamma$ acts discontinuously. The {\em limit set} $L_\Gamma$ of
$\Gamma$ is $S^2_\infty-\Omega_\Gamma$. We define
$d(N)$ to be the Hausdorff dimension of the limit set $L_\Gamma$ of $\Gamma$.

Define $\lambda_0(N)$ to be
the largest value of $\lambda$ for which there exists a positive
$C^\infty$ function $f$ on $N$ such that $\Delta f+\lambda f=0$.
(Here $\Delta=\operatorname{div}\circ \operatorname{grad}$ denotes
the Laplacian. See Sullivan 
\cite{sullivan-aspects} for a discussion of why this is equivalent
to other definitions of $\lambda_0(N)$.) Note that $\lambda_0 \ge 0$.

It is clear that if $\til N$ covers $N$ then $\lambda_0(N)\le
\lambda_0(\til N)$. 
In particular, 
$$\lambda_0(N)\le \lambda_0(\Hyp^3)=1$$
for any hyperbolic 3-manifold $N$.

If $M$ is a compact, orientable, hyperbolizable 3-manifold then we let
$TT(M)$ denote the set of complete hyperbolic 3-manifolds homeomorphic
to the interior of $M$ ($TT$ stands for ``topologically tame''; see below).
Our invariants can now be written as:
$$\Lambda(M)= \sup \{ \lambda_0(N)\mid N\in TT(M)\}$$
and 
$$D(M)=\inf \{ d(N)\mid N\in TT(M)\}.$$

We will say that an element $N =\Hyp^3/\Gamma$ of $TT(M)$
is {\em geometrically finite} if there exists  a finite collection $P$ of
incompressible annuli and tori in $\partial M$ such that 
$\hat N=(\Hyp^3\cup \Omega_\Gamma)/\Gamma$ is homeomorphic to $M-P$.

We say that a hyperbolic 3-manifold is {\em topologically tame}
if it is homeomorphic to the interior of a compact 3-manifold.
This holds by definition for each hyperbolic 3-manifold in $TT(M)$
and we remark that it is conjectured that every hyperbolic
3-manifold with finitely  generated fundamental group is topologically
tame.

\subsection{Basic properties of $\Lambda$ and $D$}
Remarkably, $\lambda_0$ and the Hausdorff dimension of the limit set
are intimately related. The following
theorem, which records that relationship, is due to Sullivan \cite{sullivan84}
in the case that $N$ is geometrically finite. If $N$ is topologically
tame but not geometrically finite, then Canary \cite{CanaryJDG} proved that
$\lambda_0(N)=0$, while Bishop and Jones 
\cite{bishop-jones} proved that if $N$ has finitely generated
fundamental group and is not geometrically finite then $d(N)=2$.

\begin{theorem}{hd-and-lambda0}{}{}
If $N$ is a topologically tame hyperbolic 3-manifold,
then $\lambda_0(N)=d(N)(2-d(N))$ unless $d(N)<1$, in which
case $\lambda_0(N)=1$.
\end{theorem}

Thus $\lambda_0$ completely determines $d$ unless $d< 1$.
The situation when $d\le 1$ was analyzed by
Sullivan \cite{sullivan79}  and Braam \cite{braam} in the case where
$N$ is convex cocompact, by Canary and Taylor \cite{canary-taylor} when
$N$ is geometrically finite,
and for all hyperbolic 3-manifolds with finitely generated fundamental
groups by Bishop and Jones \cite{bishop-jones}:

\begin{theorem}{hdone}{}
Let $M$ be a compact, orientable, hyperbolizable 3-manifold. If $N\in
TT(M)$ and $d(N)<1$, then $M$ is a handlebody or a thickened torus. If
$d(N)=1$, then $M$ is either a handlebody or an I-bundle.
\end{theorem}

The above two theorems assure
us that $D(M)$ and $\Lambda(M)$ are essentially the same invariant:

\begin{corollary}{invariants-equal}{}{}
If $M$ is a compact, orientable, hyperbolizable 3-manifold,
then $\Lambda(M)=D(M)(2-D(M))$ unless $M$ is a handlebody or a thickened
torus.
\end{corollary}

Theorem \ref{hdone} also has the following consequence:

\begin{corollary}{not-a-handlebody}{}{}
If $M$ is a compact, orientable, hyperbolizable 3-manifold which
is not a handlebody or a thickened torus, then $D(M)\ge 1$.
\end{corollary}

We note that it follows from work of Beardon \cite{beardon} that
if $M$ is a handlebody, then $D(M)=0$ and therefore $\Lambda(M)=1$.
If $M$ is  a thickened torus and $N=\Hyp^3/\Gamma\in TT(M)$, then
$L_\Gamma$ is a single point, so again $D(M)=0$ and $\Lambda(M)=1$.
This completes the first row of the scorecard.

The next proposition completes the second row:

\begin{proposition}{whenDis2}{}{}
Let $M$ be a compact, orientable, hyperbolizable 3-manifold which is not
a solid torus or a thickened torus.
Then $D(M)=2$ if and  only if any boundary component of $M$ is a torus.
\end{proposition}

\begin{pf}{}
We first suppose that any boundary component of $M$ is toroidal.
This implies that if $N\in TT(M)$ then $N=\Hyp^3/\Gamma$ has finite
volume (see Proposition D.3.18 in \cite{benedetti-petronio}).
In this case $L_\Gamma=S^2_\infty$, which implies that $d(N)=2$ and hence
$D(M)=2$.

We now suppose that $M$ has a non-toroidal boundary component.
In this case, $N$ has infinite volume for any $N\in TT(M)$.
Thurston's geometrization theorem (see \cite{Morgan}) guarantees
there is a geometrically finite manifold $N$ in $TT(M)$.
Sullivan \cite{sullivan84} and Tukia \cite{Tukia} proved that if
$N$ is geometrically finite and has infinite volume,
then $d(N)<2$. Hence $D(M)<2$.
\end{pf}

It is also  useful to note that  $D(M)$ and $\Lambda(M)$ behave monotonically
under passage to covers.

\begin{proposition}{covers}
Let $M$ and $M'$ be compact, orientable, hyperbolizable 3-manifolds, such
that the interior of $M'$ covers the interior of $M$. Then $D(M)\ge D(M')$
and $\Lambda(M)\le\Lambda(M').$ 
\end{proposition}

\begin{pf} If $N=\Hyp^3/\Gamma$ is any hyperbolic 3-manifold homeomorphic
to the interior of $M$, then it has a cover $N'=\Hyp^3/\Gamma'$ which is
homeomorphic to the interior of $N'$. Since $\Gamma'\subset \Gamma$,
$L_{\Gamma'} \subset L_\Gamma$, so $d(N)\ge d(N')$.
The assertion that $D(M)\ge D(M')$ then follows immediately from the
definition of  our invariant $D$.

The proof of the assertion that $\Lambda(M)\le\Lambda(M')$  is similar.
\end{pf}

\subsection{Outline of proof of the main theorems}
We will break the argument up into several steps. 
First we reduce to considering manifolds with incompressible
boundary. 

We say that $M$ is obtained from two manifolds $M_0$ and $M_1$ by
{\em adding a 1-handle} if $M$ is obtained from $M_0$, $M_1$ and
$D^2\times [0,1]$ 
by identifying $D^2\times \{ i\}$ with an embedded disk in $\partial M_i$
(for $i=0,1$.) $M$ is said to be obtained from $M_0$ by adding a 1-handle
if $M$ is obtained from $M_0$ and $D^2\times [0,1]$ by identifying
$D^2\times \{ 0\}$ and $D^2\times \{ 1\}$ with disjoint embedded disks
in $\partial M_0$. $M$ is said to be obtained from $\{ M_1,\ldots,M_n\}$
by adding 1-handles if it is obtained by applying the above two topological
operations finitely many times using the manifolds
$\{ M_1,\ldots,M_n\}$ as building blocks.

Bonahon \cite{bonahon83} and McCullough-Miller \cite{mccullough-miller}
showed that if $M$ is a compact irreducible 3-manifold,
then there exists a collection
$\{ M_1,\ldots, M_n\}$ of submanifolds of $M$ such
that $M$ is obtained from $\{M_1, \ldots,M_n\}$ by adding 1-handles
and each $M_i$ has incompressible boundary. 
(The boundary of a 3-manifold $M$ is incompressible if the fundamental
group of any component injects in $\pi_1(M)$).
The union $\union M_i$
is called the {\em incompressible core} of $M$.
If $M$ is a handlebody its incompressible core is a
ball, and otherwise we will assume that no component of the incompressible
core is a ball. With this convention, the incompressible core
is unique up to isotopy.

In section \ref{compress} we will use work of Patterson
\cite{Patterson-Klein} to show that one may analyze
our invariants on $M$ simply by studying the invariants
of the components of the incompressible core of $M$:

\begin{theorem}{hard-core}{}{}
Let $M$ be a compact, orientable, hyperbolizable 3-manifold.
If $\{ M_1,\ldots,M_n\}$ are the components of the incompressible core of $M$,
then $$\Lambda(M)=\min\{\Lambda(M_1),\ldots,\Lambda(M_n)\}.$$
\end{theorem}

The idea is to ``pull apart'' the groups uniformizing the components
of the incompressible core.  For example, if $M$ is obtained from $M_0$ and $M_1$
by adding a 1-handle, we first find $N_0 \in TT(M_0)$ and $N_1\in TT(M_1)$ with
$\lambda_0(N_i)$ near $\Lambda(M_i)$. We then construct a sequence $\{ N^i\}$ of
hyperbolic 3-manifolds homeomorphic to the interior of $M$
by removing half-spaces from $N_0$ and
$N_1$ which lie farther and farther from a fixed point in each, and
gluing the boundary planes together. Patterson's result is used to
show that $\{\lambda_0(N^i)\}$ converges to
$\min\{\lambda_0(N_0),\lambda_0(N_1)\}$ and hence that
$\Lambda(M)=\min\{\Lambda(M_0),\Lambda(M_1)\}$.

We now turn to 3-manifolds
with incompressible boundary, for which we will need a bit more notation.
If $X$ is an annulus or torus, we say that a map
$f:(X,\partial X)\to (M,\partial M)$ 
is {\em essential} if $f_*:\pi_1(X)\to \pi_1(M)$ is injective and
$f$ is not properly homotopic to a map of $X$ into $M$
with image in $\partial M$.
We will say that an embedding $f:R\to M$ of an $I$-bundle $R$ into $M$
is {\em admissible} if $f^{-1}(\partial M)$ is the associated
$\partial I$-bundle of $R$.
We will say that an embedding  $f:R\to M$
of a Seifert-fibred space $R$ into $M$ is admissible if $f^{-1}(\partial M)$ 
is a collection of fibres in $\partial S$.
In either case, we will say that $f$ is {\em essential}
if $f_*:\pi_1(R)\to \pi_1(M)$ is injective and whenever $X$ is a component of
$\partial R-f^{-1}(\partial M)$ then $f|_X$ is  an essential map
of a torus or annulus into $M$.

A compact submanifold $\Sigma$ of $M$  is said to be a
{\em characteristic submanifold} if
$\Sigma$ consists of a minimal collection of admissibly embedded essential
$I$-bundles and Seifert fibre spaces with the property that
every essential, admissible embedding 
$f:R\to M$ of a Seifert fibre space or $I$-bundle into $M$ is properly
homotopic to an admissible map with image in $\Sigma$.
Jaco-Shalen \cite{JS} and Johannson \cite{Johannson}
showed that every compact, orientable, irreducible 3-manifold
with incompressible boundary contains a characteristic submanifold
and that any two characteristic submanifolds are isotopic.
Hence, we often speak of {\em the} characteristic submanifold $\Sigma (M)$ of
$M$. If $M$ is hyperbolizable then every Seifert fibred component of
$\Sigma(M)$ is homeomorphic to either a solid torus or a thickened torus
(see Morgan \cite{Morgan}).

We recall that a compact irreducible
3-manifold $M$ with incompressible boundary is called a {\em generalized
book of $I$-bundles} if one may find a disjoint collection $A$ of
essential annuli in $M$ such that each component $R$ of the manifold
obtained by cutting $M$ along $A$ is either a solid torus, a thickened
torus, or homeomorphic to an $I$-bundle such that $\partial R\cap
\partial M$ is the associated $\partial I$-bundle.
One may check that a compact, orientable, irreducible 3-manifold
with incompressible boundary is a  generalized
book of $I$-bundles if the closure of any  component of
$M-\Sigma(M)$ is homeomorphic to a solid torus or a thickened torus,
and every component of $\Sigma(M)$ is a solid torus, a thickened torus,
or an $I$-bundle (see Proposition 4.3 in \cite{culler-shalen} for the
corresponding characterization of books of $I$-bundles.)

In section \ref{books} we prove:

\begin{theorem}{bookcase} 
If $M$ is a hyperbolizable generalized
book of $I$-bundles, then $\Lambda(M)=1$.
\end{theorem} 

The proof is by explicit construction. We build a hyperbolic structure
for $M$ by piecing together structures on each of the $I$-bundle
pieces, and show that if the parameters are chosen appropriately
(essentially ``pulling apart'' the $I$-bundles) the Hausdorff
dimension $d$ can be made arbitrarily close to 1. The result then
follows via the connection between $D$ and $\Lambda$.

In section \ref{unbook} we prove:

\begin{theorem}{notbook}
If $M$ is a compact, orientable, hyperbolizable 3-manifold with incompressible
boundary which is not a generalized book of $I$-bundles, then $\Lambda(M)<1$.
\end{theorem}

Here is an outline of the proof in the case that $M$ is acylindrical. 
Arguing by contradiction, we assume the existence of 
a sequence $N_i\in TT(M)$ with $\lambda_0(N_i) \to 1$.
By a theorem of
Thurston the deformation space of $M$ is compact and we can extract a
limit manifold $N$ which is homeomorphic to the interior of $M$
with $\lambda_0(N) = 1$. This
contradicts the assumption that $M$ is not a handlebody, thickened
torus or $I$-bundle, by Theorem \ref{hdone}.

If $M$ is not acylindrical, we must apply the Jaco-Shalen-Johannson
characteristic submanifold theory, and Thurston's relative compactness
theorem. In general, the limit manifold we obtain will be homeomorphic
to a submanifold of $M$, rather than to $M$ itself.

Theorems \ref{hard-core}, \ref{bookcase} and \ref{notbook} combine to
give a complete proof of Main Theorem I. Main Theorem II follows
immediately from Main Theorem I and corollary \ref{invariants-equal}.

In section \ref{remarks} we will make further comments and conjectures
about the invariants and related quantities.

\section{Reduction to the incompressible case}
\label{compress}

In this section we prove Theorem \ref{hard-core}, which  assures us that our
invariants are determined by their value on the components of the incompressible
core. We do this by showing that if $M$ is obtained by adding a 1-handle
to $M_0$ and $M_1$ (or by adding a 1-handle to $M_0$)
then $\Lambda(M)=\min\{ \Lambda(M_0),\Lambda(M_1)\}$ (or $\Lambda(M)=\Lambda(M_0)$).

The topological operation of adding a 1-handle is realized geometrically
by Klein combination.
Throughout this section we work in the 
ball model of ${\Hyp}^{3}$; the Euclidean boundary in this model
is $S^{2}$, and $\overline{{\Hyp}^{3}} ={\Hyp}^{3} \cup  S^{2}$.
If $F \subset  {\Hyp}^{3} $ is a  (convex) {\em fundamental  polyhedron}
of $\Gamma$, then $int(F)$ denotes the interior of $F$, 
$\overline{F}$ denotes the Euclidean closure of $F$, and  $F^c=\Hyp^3-F$.  
We refer the reader to sections IV.F and VI.A of \cite{maskit} 
for a full discussion of fundamental polyhedra.

\begin{theorem}{Klein-combination}{(Klein Combination, Theorem VII.A.13 in
\cite{maskit})}{} 
Let  $\Gamma_{0}$ and  $\Gamma_{1}$ be discrete subgroups
of $\PSL 2 (\C)$.  Suppose there are (convex) fundamental polyhedra  $F_{i}$
for $\Gamma_{i}$  $(i=0,1)$, so that   $int(F_0) \cup  int(F_1) =  \Hyp^{3}$.
Then the group $\Gamma$ generated by $\Gamma_0$ and $\Gamma_1$ is discrete and
is isomorphic  to $\Gamma_0*\Gamma_1$. Moreover, $F= F_0 \cap  F_1$
is a fundamental polyhedron for $\Gamma$.
\end{theorem}

If $\Gamma_0$ and $\Gamma_1$ satisfy the hypotheses of Theorem
\ref{Klein-combination} we will say that they are {\em  Klein-combinable}.
 
The main tool in the proof of Theorem \ref{hard-core}
is a result of Patterson. Let $\Gamma_{0}$ and $\Gamma_{1}$ be two
Klein-combinable  groups and suppose $N_{i} = \Hyp^3/\Gamma_{i}.$
The intuitive content of 
Patterson's result is that, if we 
``pull $\Gamma_1$ away from $\Gamma_0$'' by a suitable
sequence of conjugations $h_k\Gamma_1 h_k^{-1}$, then $\lambda_0$ of
the quotient of the combination of $\Gamma_0$ and $h_k\Gamma_1 h_k^{-1}$ 
approaches $\min\{\lambda_{0}(N_{0}),\lambda_{0}(N_{1})\}$.

The statement we give is a
version of Theorem 1  in
\cite{Patterson-Klein}. 
(Patterson actually proves his result for the critical exponents of the
groups involved. However, see \cite{sullivan-aspects}, the critical
exponent of $\Gamma$ determines $\lambda_0(\Hyp^3/\Gamma)$ and we have
translated Patterson's result into a result about $\lambda_0$.)
Let $ |g^{'}(x)|$  denote the Euclidean norm of the derivative of $g$ at $x$,
where $g \in \PSL 2(\C)$  and $x \in \Hyp^3$. 
Let $d(A_{1},A_{2})$ denote the Euclidean distance between sets
$A_{1}$ and $A_{2}$ in  $\overline{\Hyp^3}$.   

\begin{theorem}{patterson1}{(Patterson)}{}
Let $\Gamma_{0}$ and $\Gamma_{1}$ be discrete,
torsion-free subgroups of $\PSL 2(\C)$ with convex fundamental polyhedra $F_0$
and $F_1$.  Let $N_i=\Hyp^3/\Gamma_i$. Suppose there is a sequence $\{ h_k\}_{k\in \Z_+}$ in
$\PSL 2 (\C)$  so that $int(F_0)\cup h_k(int(F_1))=\Hyp^3$ for all $k$ and
$$ \frac{\sup_{w \in F_{1}^{c}} |h_{k}^{'}(w)|}{d(F_{0}^{c}, h_{k}F_{1}^{c})} \rightarrow 0 $$  
as $k \rightarrow \infty$. 
Let $\Gamma^k$ be the discrete group generated by $\Gamma_0$ and
$h_{k} \Gamma_{1} h_{k}^{-1}$ and let $N^k=\Hyp^3/\Gamma^k$. Then
$$ \lim_{k\rightarrow \infty} \lambda_0(N^k)= \min \{\lambda_0(N_0),\lambda_0(N_1)\}.$$
\end{theorem}

We begin by studying the case where $M$ is obtained from $M_0$ and $M_1$
by adding a 1-handle.
Combining Theorems \ref{Klein-combination} and \ref{patterson1}, we obtain:

\begin{proposition}{amalg-case}{}{}
Let $M_0$ and $M_1$ be two compact, orientable,
hyperbolizable 3-manifolds. If $M$ is
hyperbolizable and is obtained from $M_0$ and $M_1$ by adding a 1-handle then
$\Lambda(M) = \min \{\Lambda(M_{0}), \Lambda(M_{1})\}$.
\end{proposition}

\begin{pf}{}
We first note that Proposition \ref{covers}
implies that 
$$\Lambda(M)\le \min\{\Lambda(M_0),\Lambda(M_1)\},$$
since the interior of $M$ is covered by both the interior of $M_0$ and
the interior of $M_1$. 

Without  loss of generality assume  $\Lambda(M_{0}) =
\min\{\Lambda(M_0),\Lambda(M_1)\}$. 
Then  Proposition 2.5 ensures that $\Lambda(M_0) >0.$ 
Fix a positive $\epsilon < \Lambda(M_0)$, 
and choose $N_i = \Hyp^3/\Gamma_i \in TT(M_i)$ so  that
$\lambda_0(N_1) \ge \lambda_0(N_0) >\Lambda(M_{0}) - \epsilon >0$.
Canary \cite{CanaryJDG} showed that if a complete hyperbolic 3-manifold 
is topologically tame but not geometrically
finite, then $\lambda_0=0$.
Therefore, $N_0$ and $N_1$ are both geometrically finite.

It follows that
there exist homeomorphisms $\psi_i:\hat N_i\to M_i-P_i$, 
where $\hat N_i= (\Hyp^3\cup \Omega_{\Gamma_i})/\Gamma_i$ and $P_i$ is a
collection of disjoint incompressible annuli and tori in $\partial M_i$.
Because $M$ is hyperbolizable,  the attaching disks of the 1-handle are not
contained in toroidal boundary components of $M_i$. Hence, we may choose
disks $D_i$ in $M_i-P_i$,  such that $M$ is obtained from $M_0$ and
$M_1$ by attaching a 1-handle to the disks $D_0$ and $D_1$.
We can also assume that the pre-images $\psi_i^{-1}(D_i)$ are ``round''
disks (i.e. they are the quotients of round disks in $\Omega_{\Gamma_i}$.)

We choose lifts $\til D_i$ of $\psi^{-1}(D_i)$ to $\Omega_{\Gamma_i}$.
We may assume, by conjugating $\Gamma_1$, that the round disks $\til D_0$ and
$\til D_1$ intersect only along their common boundary circle $J$.
One can find (convex) fundamental polyhedra $F_i$ for $\Gamma_i$
such that $\til D_i$ are contained in the interiors of the intersection
of the Euclidean closures of $F_i$ with $S^2.$
Therefore $ H_i\subset int(F_i)$,  where $H_i$ denotes
the closed half-space for which $\overline{H_{i}} \cap S^2 = \til D_i$.
Thus, $\Gamma_0$ and $\Gamma_1$ are Klein-combinable.
Since $F_0\cap F_1$ is a fundamental polyhedron for  the group
$\Gamma$ generated by $\Gamma_0$ and $\Gamma_1$,  we see that
$N=\Hyp^3/\Gamma\in TT(M)$. 

Fix points $z_i$ in the interior of $\til D_i$.
Let $\gamma$ be a hyperbolic  M\"{o}bius transformation  
with  $z_0$ as its attracting fixed point  and $z_1$  as its repelling
fixed point.  We may further choose $\gamma$ so that $\gamma(\til D_0)$
is contained in the interior of $\til D_0$.
We will apply Theorem \ref{patterson1} to
the sequence $\{ h_k=\gamma^k\}$.

Note that $\gamma^k(\til D_0) \subset\gamma(\til D_0) \subset int \til D_0$
for all $k\ge 1$ which implies that
$\gamma^k(H_0) \subset \gamma(H_0)\subset H_0$.
Since $\gamma^k(F_1)$ is a fundamental polyhedron for
$\gamma^k\Gamma_1\gamma^{-k}$ and $\gamma^k(F_1)^c$ is contained
in $\gamma^k(H_0) \subset int(H_0)\subset F_0$, we see that $\Gamma_0$ and
$\gamma^k \Gamma_1 \gamma^{-k}$ are Klein-combinable and that
$F_0\cap \gamma^k(F_1)$ is a fundamental polyhedron for the
group $\Gamma^{k}$ generated by $\Gamma_{0}$ and
$\gamma^{k}\Gamma_{1}\gamma^{-k}$.
In particular, one may readily observe that $N^k=\Hyp^3/\Gamma_k\in TT(M)$
for all $k$. 

Since $F_0^c\subset H_1$ and $\gamma^k(F_1)^c\subset \gamma^k(H_0)\subset
\gamma(H_0)\subset int(H_0)$,
we see that $d(F_0^c,\gamma^k(F_1^c))\ge \delta$ for all $k\ge 1$,  where 
$\delta = d(H_1,\gamma(H_0))>0$.
 
One may  easily check that $\{(\gamma^{k})^{'}\}$ converges to $0$ uniformly on
all compact subsets of $\overline{\Hyp^3}-\{ z_1\}$.
Since $F_1^c\subset H_0$ and $\overline{H_0}$ is a compact subset of
$\overline{\Hyp^3}-\{ z_1\}$, then $\sup_{w \in F_{1}^{c}} |(\gamma^{k})^{'}(w)|$ converges to 0.

We may combine the observations above to establish that
$$\frac{\sup_{w \in F_{1}^{c}} |(\gamma^{k})^{'}(w)|}{d(F_{0}^{c}, h_{k}F_{1}^{c})} \rightarrow 0.$$
Theorem \ref{patterson1} then allows us to conclude that
$$\lim_{k \rightarrow \infty} \lambda_0(N^k)  =  \min \{\lambda_0(N_0),\lambda_0(N_1)\} \ge \Lambda(M_0)-\epsilon.$$
Therefore, $\Lambda(M)\ge\min \{\Lambda(M_0),\Lambda(M_1)\}-\epsilon.$
Since $\epsilon$ can be 
arbitrarily small, this completes the proof of Proposition \ref{amalg-case}.
\end{pf}

We will also need the following direct analogue of Theorem 1 in
\cite{Patterson-Klein}, which can be deduced from Patterson's
arguments.  

\begin{theorem}{patterson2}{(Patterson)}{}
Let $\Gamma_{0}$  be a torsion-free discrete subgroup of
$\PSL 2(\C)$ with convex fundamental polyhedron $F_0$,  and let
$\{ h_{k}\}$ be an infinite sequence in $\PSL 2 (\C)$. 
Suppose that
$F_k$ is a convex fundamental polyhedron for $\langle h_k\rangle$ such that
$int(F_0)\cup int( F_k)=\Hyp^3$ for all $k$, and 
there exists a $\delta > 0$ so that
$d(F_{0}^{c}, F_{k}^{c}) \geq \delta.$  Also, assume there
exists a $w$ lying in $F_{0} \cap F_{k}$ for all index $k$, such that for any
fixed $s > 0$,
$$ \sum_{j \neq 0}| (h_{k}^{j})'(w)|^{s} \rightarrow 0 $$
as $k \rightarrow \infty.$
Denote by $\Gamma^{k}$ the discrete group generated by $\Gamma_{0}$ and 
$\langle h_{k}\rangle$, and let
$N_{0} = \Hyp^3/\Gamma_{0}$ and $N^{k} = \Hyp^3/\Gamma^{k}.$  Then
$$ \lim_{k \rightarrow \infty} \lambda_{0}(N^{k}) = \lambda_{0}(N_{0}).$$
\end{theorem}

Theorem \ref{patterson2} allows us to handle the case where $M$ is
obtained from $M_0$ by adding a 1-handle.

\begin{proposition}{hnn}{}{}
Let $M_{0}$ be  a compact, irreducible, hyperbolizable $3$-manifold
with non-empty boundary. If $M$ is a compact hyperbolizable
3-manifold obtained by adding a 1-handle to $M_0$, then
$\Lambda(M)=\Lambda(M_0)$.
\end{proposition}

\begin{pf}{}
As in Proposition 3.3, we observe that Proposition \ref{covers} implies that
$\Lambda(M)\le \Lambda(M_0)$ and Proposition 2.5 guarantees
that  $\Lambda(M_0) > 0$.

Fix  a positive $\epsilon < \Lambda(M_0)$, and
choose $N_0\in TT(M_0)$ such that $\lambda_0(N_0)\ge \Lambda(M_0)-\epsilon$. 
As before, by \cite{CanaryJDG}, $N_0$ is necessarily 
geometrically finite. Thus,
there exists a collection $P$ of disjoint incompressible annuli and
tori in $\partial M_0$ and a homeomorphism $\psi:\hat N_0 \to M_0-P$,
where $\hat N_0 = (\Hyp^3 \cup \Omega_{\Gamma_{0}})/\Gamma_0.$  Since
$M$ is hyperbolizable, the attaching disks of the 1-handle  are not
contained in toroidal boundary components of $M$.
Hence, we may choose disks $D_0$ and $D_1$ in $\partial M_0 -P$ such that
$M$ is formed by attaching a 1-handle to the disks $D_0$ and $D_1$.
Since the interior of $\bar F\cap \Omega_{\Gamma_0}$ is a fundamental domain
for the action of $\Gamma_0$ on $\Omega_{\Gamma_0}$
(Proposition VI.A.3 in \cite{maskit}), we may further choose
$D_0$ and $D_1$ so that there are lifts $\til D_0$ and $\til D_1$ of
$\psi^{-1}(D_0)$ and $\psi^{-1}(D_1)$ which are round disks
 in the interior of  $\bar F\cap \Omega_{\Gamma_0}$.

Find a loxodromic element $\gamma$ that takes the exterior of $\til D_0$
to the interior  of $\til D_1$.
Let $H_i$ be the closed half-spaces whose Euclidean closures intersect
$S^2$ in $\til D_i$.
Then the region $F_k=\Hyp^3-(H_0\cup \gamma^{k-1}(H_1))$
is a fundamental polyhedron for $\langle \gamma^k\rangle$. Since
$int(F_k)\cup int(F_0)=\Hyp^3$, 
$\Gamma_0$ and $\langle\gamma^k\rangle$ are Klein-combinable and $F_k
\cap F_0$ is a 
convex fundamental polyhedra for the group
$\Gamma^k$ generated by $\Gamma$ and $\langle\gamma^k\rangle$. It is now easy to check that
$N_k=\Hyp^3/\Gamma^k\in TT(M)$.  

We will apply Theorem  \ref{patterson2}
with $\{h_k= \gamma^{k}\}$.  Let $\delta = d(F^c_{0}, H_{0} \cup H_{1}) > 0$
(recall $H_{i} \subset F_0$).
Because the $H_{i}$ are disjoint and 
 $F_k^c \subset H_0 \cup H_1$, then 
$$d(F_{0}^{c},F_{k}^{c}) \ge \delta$$
for all $k$.


Fix $w\in F_0\cap F_k$ for all $k>0$ and fix $s>0$.
It is well-known that $\sum_{j\ne 0} |(\gamma^j)'(w)|^s$ is finite.
It follows immediately that
$\sum_{j\ne 0} |(h_k^j)'(w)|^s=\sum_{j\ne 0}|(\gamma^{jk})'(w)|^s$ converges
to 0 as $k$ converges to $\infty$.

Thus, Theorem \ref{patterson2} implies that
$$\lim_{k\rightarrow \infty}\lambda_0(N^k)=\lambda_0(N_0).$$
Therefore, $\Lambda(M)>\Lambda(M_{0})  - \epsilon.$
Since $\epsilon>0$ was chosen arbitrarily, we have completed the proof of
proposition \ref{hnn}.
\end{pf}

Notice that one need only apply  Propositions \ref{amalg-case} and \ref{hnn} 
finitely many times in order to prove Theorem \ref{hard-core}.

\newcommand{\dist}{\operatorname{dist}}
\section{Generalized books of $I$-bundles}
\label{books}

In this section we will prove Theorem \ref{bookcase}, which says that 
$\Lambda(M)=1$ for any hyperbolizable generalized book of $I$-bundles $M$.
The key step in the proof is:

\begin{theorem}{books of I-bundles}
If $M$ is a hyperbolizable generalized book of $I$-bundles then for
any $\alpha > 1$, there exists a hyperbolic manifold $N$ homeomorphic
to $int(M)$ with $d(N) < \alpha$.
\end{theorem}

\begin{pf*}{Proof of Theorem \ref{bookcase}}
Let $M$ be a hyperbolizable generalized book of $I$-bundles which is
not a thickened torus.
Theorem \ref{books of I-bundles} implies that $D(M)\le 1$.
On the other hand, since $M$ is not a handlebody or a thickened torus,
Corollary \ref{not-a-handlebody} guarantees that $D(M)\ge 1$.
Hence, $D(M)=1$ and we conclude that $\Lambda(M)=1$
by applying Corollary \ref{invariants-equal}.
If $M$ is a thickened torus we have already observed that $\Lambda(M)=1.$
\end{pf*}

The remainder of the section is taken up with the proof of
Theorem \ref{books of I-bundles}.

\begin{pf*}{Proof of Theorem \ref{books of I-bundles}}
The characteristic submanifold of $M$ is a union of solid tori,
thickened tori, and $I$-bundles whose bases have negative Euler
characteristic.  For each $I$-bundle the subbundle over the boundary
of the base
surface is a union of annuli, which are glued to the boundary of a
solid torus or thickened torus (for a thickened torus, note that only
one of its boundaries participates in the gluing).
The union of the bases of the $I$-bundles, with boundaries glued
together inside each
solid torus, and torus boundaries of the thickened tori, comprise a
``spine'' for $M$, that is a
2-complex of which $M$ is a regular neighborhood.

We shall put a hyperbolic structure on the interior of $M$, for
which each $I$-bundle base determines a Fuchsian or extended Fuchsian
group, and each  thickened torus corresponds to a rank-2 parabolic group.
(We recall
that a discrete subgroup of $\PSL 2 (\C)$ is {\em Fuchsian} if it preserves
a half-space in $\Hyp^3$ and {\em extended Fuchsian} if it has a Fuchsian
subgroup of index 2. In either case there is a totally geodesic hyperplane
preserved by the group.) By changing
the parameters of this construction we will obtain Hausdorff dimension
arbitrarily close to 1. 

For each solid torus, the cores of annuli glued to it describe some
number $m$ of
parallel $(p,q)$ curves, for some relatively prime $p,q$ (where
$(1,0)$ denotes a meridian and $p$ is well-defined mod $q$). Consider the
following hyperbolic structure on this solid torus: Begin
with a geodesic $L$ in $\Hyp^3$ which is the boundary of $mq$
half-planes equally spaced around it (more generally the angles
between them can vary, but we will avoid this for ease of exposition).
Let $\gamma$ be a loxodromic with axis $L$, translation distance
$\ell/q$, for some (small) $\ell>0$, and rotation angle $2\pi p/q$.  The
quotient of a neighborhood of $L$ by $\gamma$ is a solid torus, which
the quotients of the half-planes meet in a collection of annuli with
boundaries glued together at the core. The intersection of these
annuli with the torus boundary give the $m$ desired $(p,q)$ curves.

For each thickened torus we choose (large) $d>0$ and
consider a horoball in $\Hyp^3$ with a
rank 2 parabolic group acting so that a fundamental domain on the
boundary is a rectangle with one sidelength $\mu_0$ and one sidelength
$md>0$. In the horoball we consider $m$ planes orthogonal to the
boundary, parallel to the $\mu_0$ side, and equally spaced 
(by distance $d$ along the boundary). In the quotient these give $m$
parallel cusps with boundary length $\mu_0$.  Here $\mu_0$ denotes a
fixed number less than the Margulis constant for $\Hyp^2$.

Choose a list of parameters $\{\ell_i\}$ for the solid tori and
$\{d_i\}$ for the thickened tori. For each base surface $S$, let $S'$
denote $S$ minus the boundary components that attach to thickened
tori, and choose a finite-area hyperbolic structure on $S'$ so that a
neighborhood of each missing boundary component is a cusp, and each
remaining boundary component that glues to a solid torus with
parameter $\ell_i$ is a geodesic of length $\ell_i$.
For each base surface we can find a Fuchsian or extended 
Fuchsian group such that the convex core of its quotient realizes
the given hyperbolic structure. (The convex core of a hyperbolic
manifold is the quotient of the convex hull of the limit set by the
associated group action.)
Note that the boundary components correspond to pure translations.
We then truncate each cuspidal end so that
the boundaries corresponding to thickened tori are horocycles of
length $\mu_0$.
We may obtain an incomplete structure on each $I$-bundle by considering
the embedding of our truncated region in the full quotient of the
associated Fuchsian or extended Fuchsian group.

For each solid torus we can then identify
neighborhoods of the corresponding boundaries of $I$-bundle bases to the
annuli arranged 
around its core, and for each thickened torus we can glue the
horocycles to the boundaries of the cusps embedded in the horoball. 
This extends consistently to the thickenings of the $I$-bundle bases
so that  we obtain
an (incomplete) hyperbolic structure on the interior of $M$, in which each
$I$-bundle base is totally geodesic.
With proper choice of the parameters, we will show that
this gives rise to a complete structure.

The developing map for this structure maps the universal cover $\til M$ to
$\Hyp^3$ by a locally isometric immersion (see e.g. Benedetti-Petronio
\cite{benedetti-petronio} \S B.1). Let $\Gamma$ denote the holonomy group.
Each component of the lift $\til S$ of a base surface $S$
maps to a totally geodesic subset of $\Hyp^3$. These subsets, which we
will call {\em flats}, are arranged in a ``tree'', in this sense: 
For a given flat $F$, at each lift of a geodesic boundary of its base
surface there is a collection of $mq-1$ other flats, equally spaced
(where $m$ is the number of annuli glued to the corresponding torus,
and $(p,q)$ describes the slopes of these annuli, as above).
At each parabolic fixed point corresponding to one of the boundaries
glued to a thickened torus, there is an bi-infinite sequence of other
flats, arranged with equal spacing around a horoball based at this point.
The corresponding graph of adjacencies is a tree (of infinite valence).

\subsection{Discrete holonomy.}
\label{discrete holonomy}
Let $\ell_0 = \max \ell_i$ and $d_0 = \min d_i$. Let $\theta_0 = \min
2\pi/q_im_i$ where $\{m_i\}$ and $\{(p_i,q_i)\}$ describe the gluings
for the solid tori.  We will show that, if $\ell_0$ is
sufficiently small and $d_0$ sufficiently large, the holonomy group $\Gamma$
is discrete, and the quotient manifold is homeomorphic to $M$.

We first make the following geometric observation, which is a standard
type of fact for broken geodesics in $\Hyp^n$. 

\begin{lemma}{broken geodesics}
Given $\theta\in(0,\pi]$ there exists $K\ge 0$ such that the following
holds.  Let $\gamma$ be a broken geodesic in $\Hyp^3$ composed of a
chain of $n$ segments $\gamma_1,\ldots,\gamma_n$ of lengths $k_i > K$
that meet at angles $\theta_i \ge \theta$.  Let $P_i$ denote the
orthogonal bisecting plane to $\gamma_i$. Then the $P_i$ are all
disjoint, and each $P_j$ separates $P_i$ and $\gamma_i$ from $P_k$ and
$\gamma_k$ whenever $i<j<k$.  Furthermore $\dist(P_i,P_{i+1}) \ge
\half(k_i+k_{i+1})-K.$
\end{lemma}
\begin{pf}
Choose $K$ by the formula
$$
\cosh^2 K/2 = {2\over 1-\cos \theta}.
$$
A little hyperbolic trigonometry shows that if two segments meet at
their endpoints at
angle $\theta$ then the planes orthogonal to the segments at a
distance $K/2$ from the intersection point meet at a single point at
infinity, and if the angle is greater than $\theta$ the planes are
disjoint.

Now consider for each segment of $\gamma$ the family of planes
orthogonal to it, excluding the ones closer than $K/2$ to either
endpoint (a nonempty family since $k_i>K$). The planes meeting 
any segment thus separate the planes meeting the previous segment from
those meeting the next segment,
and the distance between the first and last plane
for segment $i$ is $k_i-K$. 
The statement for the bisecting planes follows from this.
\end{pf}

Recall that the {\em $\mu$-thin part} of a flat $F$ denotes the points
where some element of the stabilizer of $F$ acts with translation
$\mu$ or less. If $\mu$ is smaller than the Margulis constant, this
set consists of a union of disjoint pieces, each of which is either a
horodisk around a parabolic fixed point or a neighborhood of an axis
of a translation. The $\mu$-thick part is the complement of the
$\mu$-thin part.

Given any two points $x,y$ in two flats $F,F'$,
let $F=F_1,\ldots,F_n=F'$ denote the sequence of flats in the tree connecting
them. Each $F_i$ and $F_{i+1}$ share either a geodesic boundary or a
parabolic fixed point at infinity, called $F_i\intersect F_{i+1}$ in
either case. There is a chain of geodesics $\{\alpha_i\}$
connecting $x$ to $y$ such that $\alpha_i\subset F_i$, and $\alpha_i$
meets $\alpha_{i+1}$ at $F_i\intersect F_{i+1}$ (possibly at infinity).
The chain is uniquely determined by the condition that, whenever
$F_i\intersect F_{i+1}$ is a geodesic boundary, $\alpha_i$ meets it
orthogonally. 
Whenever $\alpha_i,\alpha_{i+1}$ meet in a parabolic point,
adjust them as follows: Truncate each at the point where it enters the
$\mu_1$-horoball of the corresponding parabolic group ($\mu_1<\mu_0$ will
be determined shortly), and join the new endpoints with a geodesic,
which we note makes an angle greater than $\pi/2$ with $\alpha_i$ and
$\alpha_{i+1}$. When $i=n-1$, $\alpha_n$
may be entirely contained in the $\mu_1$-horoball,
and in that
case we remove $\alpha_n$ entirely and join $y$ directly to the
truncated $\alpha_{n-1}$.
Call the resulting chain of geodesics $\gamma_{x,y}$.

Suppose that $x$ is in the $\mu_0$-thick part of $F$. We claim that,
given any $k$, 
if $\ell_0$ and $\mu_1$ are sufficiently small and $d_0$ is
sufficiently large, each segment of $\gamma_{x,y}$,
except possibly the last, has length at least $k$.
By the collar lemma for hyperbolic surfaces, if $\ell_0$ and $\mu_1$ are
sufficiently small then the $\mu_0$-thick part of each quotient
surface is separated from its boundary by at least $c\log \mu_0/\ell_0$,
and from the $\mu_1$-thin parts of the cusps by $c\log \mu_0/\mu_1$, for a
fixed constant $c$. This bounds from below the length of each
segment in $\gamma_{x,y}$, except possibly the last segment containing
$y$, and the additional segments added in horoballs. Each segment of
the latter type has length at least $c\log d_0\mu_1/\mu_0$ (for a
constant $c$) since the horospherical distance between flats on the
boundary of the $\mu_0$-horoball is a multiple of $d_0$ by construction. 
Thus, choosing $d_0$ large enough (after the choice of $\mu_1$ is made) this
gives a high lower bound for the horoball segments, 
and establishes our claim.

Any two segments meet at
angle at least $\theta = \min\{\theta_0,\pi/2\}$, so let $K=K(\theta)$
be the constant given in Lemma \ref{broken geodesics} and suppose
$k\ge 2K$ and $\ell_0, d_0$ and $\mu_1$ are determined as above.
Lemma \ref{broken geodesics} then provides a sequence of planes with definite
spacing that separate $x$ from $y$. 

In particular we can deduce that
any two flats which are non-adjacent in the tree are disjoint, and
more generally,
fixing $x$ in the $\mu_0$-thick part of $F$, for any path
$F=F_1,\ldots,F_n$ of successively adjacent flats in the tree, that 
$$
\dist(x,F_n) \ge (n-2)(k-K) + k-K/2.
$$

In particular the entire tree of flats is (properly) embedded in
$\Hyp^3$, and therefore $\Gamma$ is discrete. 

It remains to show that $N=\Hyp^3/\Gamma$ is homeomorphic to $M$.
Bonahon's theorem \cite{bonahon} guarantees that $N$ is topologically
tame (we could also deduce this directly by showing that $\Gamma$ is
geometrically finite). 
Since a neighborhood of the tree of flats embeds, it must be the
homeomorphic developing image of a neighborhood of the lift to the
universal cover of the spine of $M$. It follows that $M$ embeds
in $N$, by a map which is a homotopy equivalence. 
By a theorem of McCullough-Miller-Swarup \cite{mms}, this implies that
$N$ is homeomorphic to the interior of $M$.

{\em Remark:} Another approach to this construction is by
means of  Klein-Maskit combinations (see Maskit \cite{maskit-panelled}.)

\subsection{Hausdorff dimension.}
\label{Hausdorff dimension}
We next show that, with further restrictions on  $\ell_0$ and $d_0$,
we can obtain upper bounds on Hausdorff dimension. This will be done
directly, by exhibiting an appropriate family of coverings.

Choose one flat $F_0$ as the root of the tree. Choose a point $x_0$
in the $\mu_0$-thick part of $F_0$.
Normalize the picture in the upper half-space model so that
the plane $H_0$ containing $F_0$ is a hemisphere meeting the complex
plane in the unit 
circle $C_0$, and so that $x_0$ is the point $(0,0,1)$.
Each child $F'$ of $F_0$ in the tree structure is of one of two types:
Type (1): if $F'$ meets $F_0$ along a geodesic $L$, $F'$ is contained in a
half-plane $H'$ which meets the complex plane in an arc $C'$ of a circle.
There 
are a finite number (at most $2\pi/\theta_0$) of other flats adjoined at $L$.
Type (2): If $F'$ meets $F_0$ along a parabolic fixed point, it is part
of a sequence of flats $\{F_n\}_{n\in\Z}$ meeting at that point, where
$F_n$ are all children of $F_0$ for $n\ne 0$. These
meet the plane in a family of concentric circles $\{C_n\}$ tangent at
the same point. These divide into $C_0$ itself, and the circles
outside and inside. Call the set of outer ones (and of inner ones) an
``earring''. 

The same description holds for the children of any flat. We thus get a
family of circular arcs $\{C\}$ arranged in a tree structure with root
$C_0$. (Note that for type (1) flats the arc is just the portion of a
circle where a half-plane meets $\C$, whereas for type (2) it is a
whole circle.)
For any $C$ let $s(C)$ denote the set of its children, and
similarly $s^2(C) = s(s(C))$, etc. Let $r(C)$ denote the diameter of
$C$, which of course is uniformly comparable to the length of $C$.

The limit set $L_\Gamma$ of $\Gamma$ is contained in the closure
$\hat L_\Gamma$ of the union of these arcs $C$.

Fix any positive $\rho< 1/2$.
We claim that we can choose the parameter
$d_0$ sufficiently large and $\ell_0$ sufficiently small, so that the
following holds (where $c_0$ is a fixed constant):

\begin{itemize}
\item For any $C$ and $D\in s(C)$, 
\begin{equation}
\label{geometric bound}         
                r(D) \le \rho r(C).
\end{equation}
\item If $D_1,D_2,\ldots$ are the nested circles in an earring, with $D_1$
the outermost, then 
\begin{equation}
\label{harmonic bound}
                r(D_n) \le c_0{r(D_1)\over n}.
\end{equation}
\end{itemize}

\begin{pf*}{Proof of (\ref{geometric bound})}
Let $\mu_2<\mu_0$ be a constant to be chosen later. 
Let $H$ and $H'$ be the hemispheres containing $C$ and $D$,
respectively. They meet either in a geodesic $g$ or at a parabolic
fixed point. Let $J(H,H')$ denote the component of the $\mu_2$-thin
part associated to the intersection in either case.
Let $x$ be the top point of $H$ (in the upper
half-space). Suppose that $x$ is outside $J(H,H')$. Then the geodesic
chain $\gamma_{x,y}$ for any $y\in H'$, constructed as in
\S\ref{discrete holonomy} but with $\mu_2$ taking the place of
$\mu_0$, has initial segment $\gamma_1$ of length at 
least $k$, where $k$ can be made arbitrarily large by making 
$\ell_0/\mu_2$ and $\mu_1/\mu_2$ small, and $d_0$ large.
Given these choices, we conclude via Lemma \ref{broken geodesics}
that all of $H'$ is separated from $x$ by the bisecting
hemisphere of $\gamma_1$, which is distance at least $k/2$ from $x$,
and is thus of Euclidean diameter at most $ce^{-k/2}\diam(H)$ for a
fixed $c$.  This gives the desired bound on $r(D)$, if $k$ is chosen
so that $ce^{-k/2}\le \rho$.

It remains to show that, with appropriate choice of $\mu_2$,  the top
$x$ of each $H$ is outside $J(H,H')$ 
for any child $H'$ of $H$. For the root of the tree this holds by our
normalization. We argue by induction. Suppose that $H_1,H_2,H_3$ are
hemispheres such that $H_i$ is the parent of $H_{i+1}$ and
$x_i$ are the tops of $H_i$. If the inductive hypothesis holds for
$H_1$ then, in particular, $\diam(H_2) \le \diam(H_1)$ by the above
paragraph. It follows that $\dist(x_2,H_1)$ is bounded from above
by a fixed number
$a$. However, if $x_2$ were contained in $J(H_2,H_3)$ then there would
be a separation between $x_2$ and $J(H_1,H_2)$ of at least
$c\log(\mu_0/\mu_2)$, since $J(H_1,H_2)$ and $J(H_2,H_3)$ meet $H_2$
in two distinct, and hence disjoint, components of the $\mu_0$ thin part.
Assuming that $\mu_2$ is sufficiently short
this distance is long enough that we can construct a geodesic chain
$\gamma_{x_2,y}$ from $x_2$ to any $y\in H_1$ with long initial segment, and
apply Lemma \ref{broken geodesics} to conclude $\dist(x_2,H_1) > a$, a
contradiction.
Thus there is an a-priori choice of $\mu_2$ which guarantees that
$x_2$ will be outside $J(H_2,H_3)$, and we are done by induction.
\end{pf*}

\begin{pf*}{Proof of (\ref{harmonic bound})} We observe that since, by
(\ref{geometric bound}), $D_1$ is at most half the size of the parent $C$, we
may re-normalize, by a M\"obius transformation whose derivative is within
a universally bounded ratio of a constant on all of the
$D_i$, so that $C$ becomes a straight line meeting the $D_i$ at the origin.
The $D_i$, and $C$, are then taken
by the map $z\mapsto 1/z$ to a sequence of equally spaced parallel lines.
An easy computation gives (\ref{harmonic
bound}), where the constant $c_0$ comes from the initial re-normalization.
\end{pf*}

\medskip

We now claim the following,
for any $2\ge \alpha>1$:
\begin{equation}
\label{power sum bound}
\sum_{D\in s(C)} r(D)^\alpha \le a_0\rho^{\alpha-1}{\zeta(\alpha)\over
2^{\alpha-1} -1}  r(C)^\alpha
\end{equation}
where $a_0$ is a fixed constant 
and $\zeta(\alpha) = \sum 1/n^\alpha$ is the usual Zeta function.

To prove this consider first the children of type (1): these are
arranged in groups of bounded number which subtend a common interval
on $C$, and any two such intervals are disjoint. The lower bound
$\theta_0$ on the angle at which any such child meets $C$ implies that
its diameter is comparable to the diameter of the interval.
The sum of
lengths of intervals is at most the length of $C$, so there is some
constant $a_1$ such that $\sum r(D) \le
a_1 r(C)$, for $D$ of type (1). Since also $r(D)\le \rho r(C)$, we 
make the following observation:
If $\sum x_i \le ax$ and each $x_i \le \rho x$, then
$\sum x_i^\alpha \le \sum x_i(\rho x)^{\alpha-1} \le a\rho^{\alpha-1}
x^\alpha$. Thus we can 
bound the contribution of type (1) children by:
\begin{equation}
\label{type 1 bound}
\sum_{\text{type 1}} r(D)^\alpha \le a_1\rho^{\alpha-1}r(C)^\alpha.
\end{equation}

For the children of type (2) the sum of lengths is infinite and we
must take more care. Consider first an earring $D_1,D_2,\ldots$ with
$D_1$ outermost. By (\ref{harmonic bound}), we have
\begin{equation}
\label{earring bound}
\sum_{n=1}^\infty r(D_n)^\alpha \le c_0^\alpha\zeta(\alpha) r(D_1)^\alpha.
\end{equation}
At each parabolic point $p$ on $C$ there are two earrings (inside and
outside the circle). Let $D_p$ denote the outermost circle of the
outside earring. Clearly it just remains to bound $\sum_p
r(D_p)^\alpha$. 

Note first that all the $D_p$ are disjoint, by the argument of 
\S\ref{discrete holonomy} showing that the tree of flats embeds.
It follows, we claim, for any $\delta>0$, that
\begin{equation}
\label{delta interval bound}
        \sum_{{r(D_p)\over r(C)}\in[\delta/2,\delta]} r(D_p)^\alpha
                        \le a_2  \delta^{\alpha-1} r(C)^\alpha.
\end{equation}
Each $D_p$ projects radially to an interval on $C$, and the condition
that diameters lie in $[r(C)\delta/2,r(C)\delta]$, together with
disjointness of the $D_p$, means that these
intervals cover 
$C$ with multiplicity at most $2$. This implies
$\sum r(D_p) \le a_2 r(C)$ for this subset with some constant $a_2$,
and (\ref{delta interval bound}) follows, using the same
observation as for the type (1) children.

Summing over $\delta = \rho/2^k$
for $k=0,1,\ldots$, we obtain a bound for the sum over all outer
circles of (outside) earrings: 
\begin{align}
\label{outer circle sum}
\sum_{p} r(D_p)^\alpha & \le a_2  r(C)^\alpha \sum_{k=0}^\infty
                                \left({\rho\over 2^{k}}\right)^{\alpha-1}
\notag\\
                      & \le {2a_2 \rho^{\alpha-1}\over 2^{\alpha-1}-1} r(C)^\alpha.
\end{align}
The same argument works for the outer circles of earrings contained
inside $C$, doubling our bound. 
Combining with the bound (\ref{earring bound}) for the sum over
each earring, we obtain the inequality (\ref{power sum bound}) for the
sum over type (2) children. Now combining with (\ref{type 1 bound}) we
get (\ref{power sum bound}) over all the children of $C$ (note that in
(\ref{type 1 bound}) the factor $\zeta(\alpha)/(2^{\alpha-1}-1)$ does
not appear, but this does not matter since it has a positive lower bound when
$\alpha\in(1,2]$ ).

\medskip

We shall now define, for any $\ep_0 >0$, a covering of the
closure $\hat L_\Gamma$ of the union of arcs $C$, by balls of radius
less than or equal to $\ep_0$. 

For each arc $C$ we shall inductively assign a number $\ep(C)$ with
which to cover $C$. Let $\ep(C_0) = \ep_0$. For a child $D$ of $C$ let
$$
\ep(D) = {r(D)\over \rho r(C)} \ep(C).
$$
Note that $\ep(D) \le \ep(C)$ by (\ref{geometric bound}). Furthermore we
observe by induction that if $C\in s^j(C_0)$, the $j$-th level of the
tree, then
$$
\ep(C)  = {r(C)\over \rho^j r(C_0)} \ep_0 \le \ep_0.
$$

Recall that $\rho<1/2$. If $D$ is a child of $C$ and $T(D)$
is the union of all arcs in the subtree whose root is $D$, we
immediately have by (\ref{geometric bound}) that $T(D)$ is contained
in a ball of radius
$$
\sum_{k=0}^\infty \rho^k r(D) < 2r(D)
$$
around any point of $D$. Thus there is some  fixed $a_3$ for which
there exists a covering of $C$ by $a_3 r(C)/\ep(C)$ 
balls of radius $\ep(C)$, which also covers the
closure of any subtree descended from $C$ with root $D$, provided
$r(D) \le \ep(C)/2$.

Now given $k>0$ let $\ep_0 = \rho^k r(C_0)$ and let $U_k$ denote the covering
which is the union of these coverings for all $C$ in levels $0$
through $k$. We claim that $U_k$ is in fact a covering of all of $\hat
L_\Gamma$. 
For, if $D\in s(C) $ and $C$ is at level $k$,
$r(D) \le \rho r(C)$, and 
$\ep(C) = (r(C)/\rho^k r(C_0)) \ep_0 = r(C)$. Thus the covering for
$C$ covers the closure of $T(D)$, as above.

We shall now compute the $\alpha$-dimensional mass of this covering.
Let $M_k$ denote the sum $\sum r_i^\alpha$ over the balls of $U_k$,
where $r_i$ denotes the radius of the $i$-th ball
and let $M_k(C)$ for any arc $C$ denote the sum over just the subset of
balls covering $C$.

For each $C$ we have
$$M_k(C) = a_3 (r(C)/\ep(C)) \ep(C)^\alpha
 = a_3 r(C) \ep(C)^{\alpha-1}.$$
If $D \in s(C)$, we get
\begin{align}
M_k(D)  &= a_3 r(D) \left({r(D) \ep(C)\over \rho r(C)}\right)^{\alpha-1}
                \notag \\
        & = a_3 \left({\ep(C)\over \rho r(C)}\right)^{\alpha-1} r(D)^\alpha.
                \notag
\end{align}
Summing over all $D\in s(C)$ and using (\ref{power sum bound}),
\begin{align}
\label{mass sum bound}
\sum_{D\in s(C)} M_k(D)  & \le a_3 \left({\ep(C)\over \rho
                r(C)}\right)^{\alpha-1} a_0 \rho^{\alpha-1}
                {\zeta(\alpha) \over 2^{\alpha -1}-1} r(C)^\alpha \notag\\ 
& = a_0 {\zeta(\alpha) \over 2^{\alpha-1}-1} a_3 r(C)
                \ep(C)^{\alpha-1}\notag\\
& = A(\alpha) M_k(C).
\end{align}
Where we abbreviate
$A(\alpha) = a_0a_3\zeta(\alpha)/(2^{\alpha-1}-1)$.
Applying this inductively, we get
\begin{equation}
\sum_{D\in s^j(C_0)} M_k(D)\le A(\alpha)^j M_k(C_0) 
\label{one level sum}
\end{equation}
For any $j\le k$.
Let us assume $A(\alpha) \ge 2$ (since we may always enlarge $a_0$).
Then $\sum_{j=0}^k A(\alpha)^j \le A(\alpha)^{k+1}$, and 
summing up (\ref{one level sum}) over levels $0$ through $k$ for the
covering $U_k$, we get 
\begin{equation}
M_k \le {A(\alpha)^{k+1}} M_k(C_0).
\end{equation}
By the choice of $\ep_0 = r(C_0) \rho^k$, we have
$$
M_k(C_0) = a_3 r(C_0)^\alpha (\rho^{\alpha-1})^k.
$$
Thus, given any $\alpha>1$, we may choose $\rho$ small enough that
$\rho^{\alpha-1} < 1/A(\alpha)$, and then $\lim_{k\to\infty} M_k =
0$, and $\hat L_\Gamma$ (hence $L_\Gamma$) has zero Hausdorff measure in
dimension $\alpha$. Thus $d(\Hyp^3/\Gamma) \le \alpha$, which concludes the
proof of Theorem \ref{books of I-bundles}.
\end{pf*}

\section{More preliminaries}
\label{more}

In this section we recall more of the background which
will be needed to handle manifolds with incompressible
boundary which are not books of I-bundles.
We first recall Bonahon's theorem
about topological tameness and some basic facts about geometric convergence.
In section \ref{deform} we recall Thurston's relative compactness
theorem and prove an unmarked version of it which will be 
the key technical tool in the proof of theorem \ref{notbook}.

\subsection{Bonahon's theorem}
\label{TT}
Bonahon \cite{bonahon} proved that
a hyperbolic 3-manifold $\Hyp^3/\Gamma$ is topologically tame
if $\Gamma$ satisfies Bonahon's condition (B), which is the following:
whenever $\Gamma=A*B$ is a non-trivial free decomposition of $\Gamma$,
there exists a parabolic element of $\Gamma$ which is not
conjugate to an element of either $A$ or $B$.

Bonahon provided
the following topological interpretation of his condition (B)
(see Proposition 1.2 in \cite{bonahon}.)

\begin{lemma}{Bcondition}{}{}
Let $\Gamma$ be a torsion-free Kleinian group.
Suppose that $M$ is a compact 3-manifold, $P$ is a collection
of homotopically non-trivial annuli in $\partial M$, no two of which
are homotopic in $M$ and every component of $\partial M-P$
is incompressible in $M$. Then, if there exists an isomorphism
$\phi:\pi_1(M)\to \Gamma$ such that $\phi(g)$ is parabolic
if $g$ is conjugate to an element of $\pi_1(P)$, then
$\Gamma$ satisfies Bonahon's condition (B).
\end{lemma}

\subsection{Geometric convergence}
\label{geom}

We say that a sequence of Kleinian groups $\{ \Gamma_j\}$ converges
{\em geometrically}
to a Kleinian group $\Gamma$ if for every $\gamma\in \Gamma$, there
exists a sequence $\{\gamma_j\in\Gamma_j\}$ converging to $\gamma$
and if any accumulation point of any sequence $\{\gamma_j\in\Gamma_j\}$
is contained in $\Gamma$.
We will need the following observation,  whose proof we sketch.
For a complete argument, from a slightly different  point of view, see Taylor
\cite{taylorwhat}.

\begin{lemma}{lambda-semicontin}{}
Suppose that $\Gamma_i$ is a sequence of torsion-free Kleinian
groups converging geometrically to a torsion-free Kleinian group $\Gamma$.
Let $N_i=\Hyp^3/\Gamma_i$ and $N=\Hyp^3/\Gamma$. Then
$$\limsup \lambda_0(N_i)\le \lambda_0(N).$$
\end{lemma}

\begin{pf}{}
Let $\{N_j\}$ be a subsequence of $\{ N_i\}$ such that
$\{ \lambda_0(N_j)\}$ converges to $L=\limsup \lambda_0(N_i)$.
Let $f_j$ be a positive  $C^\infty$-function such that
$\Delta f_j+\lambda_0(N_j)f_j=0$ and let $\til f_j$ be the lift of $f_j$ to
a map $\til f_j:\Hyp^3\to \R$. We may scale $\til f_j$ so that
$\til f_j(\vec 0)=1$ where $\vec 0$ denotes the origin of $\Hyp^3$.
Yau's Harnack inequality \cite{yau-harnack} and basic elliptic theory
guarantee that there is a subsequence of $\{\til f_j\}$ which converges
to a  positive $C^\infty$-function $\til f$ such that $\Delta \til f+L\til f=0$.
Since $\til f_j$ was $\Gamma_j$-invariant and $\Gamma_j$ converges
geometrically to $\Gamma$, $\til f$ is $\Gamma$-invariant and hence descends to
a $C^\infty$-function on $N$ such that $\Delta f+Lf=0$. It follows that
$\lambda_0(N)\ge L$.
\end{pf}

\subsection{Deformation theory of Kleinian groups}
\label{deform}
In the proof of Theorem \ref{notbook} we will consider
sequences of Kleinian groups in an algebraic deformation space. 
Let us therefore introduce the following terminology.
If $G$ is a group, let ${\cal D}(G)$ denote the space of
discrete, faithful representations of $G$ into $\PSL 2 (\C)$.
Let $AH(G)={\cal D}(G)/\PSL 2 (\C)$ where $\PSL 2(\C)$ acts
by conjugation of the image.
If $c$ is a loop in $M$ represented
by the element $g\in \pi_1(M)$ and $\rho\in AH(\pi_1(M))$, then $l_\rho(c)$
is the translation length of $\rho(g)$ if $\rho(g)$ is hyperbolic and 0 if
$\rho(g)$ is parabolic.

We will use the following basic lemma which relates convergence in
${\cal D}(G)$ and geometric
convergence (see \cite{jorgensen-marden}).

\begin{lemma}{}
Suppose that $G$ is a torsion-free, non-abelian group and 
$\{\rho_j\}$ is a sequence in ${\cal D}(G)$ which converges to 
$\rho\in {\cal D}(G)$. Then there is a subsequence of $\{\rho_j(G)\}$
which converges geometrically to a torsion-free Kleinian group
$\Gamma$ which contains $\rho(G)$.
\end{lemma}

In order to state Thurston's relative compactness theorem we need 
to introduce the {\em window} $W$ of a compact hyperbolizable 3-manifold $M$
with incompressible boundary. The window
$W$ consists of the $I$-bundle components of the characteristic
submanifold $\Sigma(M)$ together
with a thickened neighborhood of every essential annulus in
$\partial \Sigma(M) -\partial M$  which is not the boundary of
an $I$-bundle component of $\Sigma(M)$.
The window is itself an $I$-bundle over a surface $w$; $w$ is called the
{\em window base}.

The following is Thurston's relative compactness theorem
for hyperbolic structures on $M$ (see \cite{ThIII} or
Morgan-Shalen \cite{MS}.)

\begin{theorem}{}
Let $M$ be a compact, orientable, hyperbolizable
3-manifold with incompressible boundary and window $W$.
If $G$ is any subgroup of $\pi_1(M)$ which
is conjugate to the fundamental group of a component
of $M-W$ whose closure is not a thickened torus,
then the image of the induced map from $AH(\pi_1(M))$ into $AH(G)$ has 
compact closure.
\end{theorem}

We wish to extend this theorem to the following ``unmarked'' version
of Thurston's theorem.

\begin{theorem}{unmarked-thurston}{}{}
Let $\{\rho_i\}$ be a sequence in $AH(\pi_1(M))$. We may
then find a subsequence $\{\rho_j\}$, a sequence of elements
$\{\phi_j\}$ of $Out(\pi_1(M))$, and a collection
$x$ of disjoint, non-parallel, homotopically non-trivial
simple closed curves in the window base $w$ such that  if $G$ is any 
subgroup of $\pi_1(M)$ which is conjugate to the fundamental group of
a component of $M-X$ whose closure is not a thickened torus
(where $X$ is the total space of the $I$-bundle over $x$),
then $\{\rho_j\circ\phi_j|_G\}$ converges in $AH(G)$.
Moreover, if $c$ is a curve in $x$, 
then $\{ l_{\rho_j\circ\phi_j}(c)\}$ converges to 0.
\end{theorem}

The idea of the proof is the following. By Thurston's theorem, the
restrictions of the representations to the complementary components of
the window have convergent subsequences, and in particular the lengths
of the window boundaries are bounded. Now for each element of the
sequence we represent each component of the window base as a pleated
surface with geodesic boundary, and observe that the hyperbolic
structures induced on these surfaces, after appropriate remarkings and
restriction to a subsequence, either converge or develop cusps. In the
latter case we cut along the curves which become cusps (these are the
family $x$), and argue that the representations restricted to the
remaining components converge up to subsequence. Finally, we reglue
along the window boundaries which did not converge to cusps, and again
after restriction to a subsequence obtain the desired convergence.

\begin{pf}{}
The following proposition states the corresponding fact for hyperbolic
surfaces. This fact is reasonably standard,  and
one may construct a proof using techniques described in Abikoff 
\cite{abikoff} and Harvey \cite{harvey}.

\begin{proposition}{unmarked-Fuchsian}{} Let $S$ be a compact
surface with boundary and let $\{\psi_i\}$ be a sequence of discrete
faithful representations
of $\pi_1(S)$ into $Isom(\Hyp^2)$. If there exists $K$ such
that $l_{\psi_i} (\partial S)\le K$ for all $i$, then there
exists a subsequence $\{\psi_k\}$ of $\{\psi_i\}$,
a collection $y$ of disjoint, homotopically distinct curves 
in $S$ and a collection of homeomorphisms $h_k:S\to S$ which
are the identity on $\partial S$ such that if $R$ is a component
of $S-y$, then $\{\psi_k\circ (h_k)_*|_{\pi_1(R)}\}$
converges in $AH(\pi_1(R))$.
Moreover, $\{ l_{\psi_k \circ (h_k)_*}(y_0)\}$ converges to 0, for any
component $y_0$ of $y$.
\end{proposition}

We begin by using Proposition \ref{unmarked-Fuchsian} to choose $x$,
a subsequence $\{\rho_k\}$ of $\{\rho_i\}$, and a
sequence of homeomorphisms $\{f_k:w\to w\}$  such that if
$R$ is a component of $w-x$, then $\{\rho_k\circ (f_k|_R)_*\}$
converges in $AH(\pi_1(R))$.

Let $w_1,\ldots,w_m$ be the components of $w$. We choose the
restriction of $x$ to each component inductively, at the same
time passing to subsequences. 
We inductively assume that we have chosen a subsequence $\{\rho_k\}$ and
the restriction of $x$ and $f_k$ to $w_1\cup\cdots w_{l-1}$ such
that if $R$ is a component of $(w_1\cup\cdots\cup w_{l-1})-x$,
then $\{\rho_k\circ (f_k|_R)_*\}$ converges in $AH(\pi_1(R))$.

If $w_l$ is an annulus or a M\"obius strip,
let $c_l$ be a core curve of $w_l$.  Thurston's relative compactness theorem
implies that $l_{\rho_i}(c_l)$ is bounded above.
We may then pass to a subsequence  of $\{\rho_k\}$, again
called $\{\rho_k\}$, such that $\{\rho_k|_{\pi_1(w_l)}\}$ converges in $AH(\Z)$.
We include $\partial w_l$ in $x$ if and only if $\{l_{\rho_k}(c_l)\}$
converges to 0.
Let  the restriction of $f_k$ to $w_l$ be the identity map.

Now suppose that $w_l$ is not an annulus or a M\"obius strip, in
which case $w_l$ has negative Euler characteristic.  We may
pass to a subsequence of $\{\rho_k\}$, again called
$\{\rho_k\}$, such that if $z$ is any boundary
component of $w_l$ then either $l_{\rho_k}(z)=0$ for all $k$ or
$l_{\rho_k}(z)\ne 0$ for any $k$.  Let $w_l'$ be obtained from $w_l$
by removing any component $z$ of $\partial w_l$ such that
$l_{\rho_k}(z)=0$ for all $k$.    One can then find, for
each $k$, a complete finite-area hyperbolic metric $\tau_k$ on $w_l'$
with geodesic boundary, and a path-wise isometry $r_k:(w_l',\tau_k)\to
N_k$ in the homotopy class determined by the restriction of $\rho_k$ to
$\pi_1(w_l)$, such that that $r_k(\partial w_l')$ is a collection of closed
geodesics.  (A map between Riemannian manifolds
is a {\em pathwise isometry} if any rectifiable path in the domain is
taken by the map to a path of equal length.
Typically one takes each $r_k$ to be a pleated surface, see \cite{CEG} or
\cite{ThI}.  In particular, note that neighborhoods of the missing
boundaries of $w_l'$ are cusps, which map into cusps in $N_l$.)

Then $(w_l',\tau_k)$ is isometric to the convex core of $\Hyp^2/\Theta_k$
for some discrete subgroup $\Theta_k$ of $Isom(\Hyp^2)$ and
there is an induced discrete faithful representation
$\psi_k:\pi_1(w_l) \to Isom(\Hyp^2)$ with image $\Theta_k$ such
that $l_{\psi_k}(c)\ge l_{\rho_k}(c)$ for any simple
closed curve $c$ in $w_l$ and  all $k$.

Thurston's relative compactness theorem guarantees that there exists $K$
such that $l_{\psi_k}(\partial w_l')=l_{\rho_k}(\partial w_l)\le K$ for all
$k$.  Then $\{\psi_k\}$ satisfies the hypotheses of Proposition
\ref{unmarked-Fuchsian}. Hence, after perhaps passing to
another subsequence of $\{\rho_k\}$, again called $\{\rho_k\}$, we
obtain a collection $y$ of simple closed curves in $w_l$ and a
sequence of homeomorphisms $h_k:w_l\to w_l$ such that if $R$ is a
component of $w_l-y$, then $\{\psi_k\circ (h_k|_R)_*\}$ converges in
$AH(\pi_1(R))$.  We then append $y$ to $x$ and we also add to $x$ any
component $z$ of $\partial w_l$ such that $\{ l_{\rho_k}(z)\}=\{
l_{\psi_k}(z)\}$ converges to 0. Let the restriction of $f_k $ to
$w_l$ agree with $h_k$.  If $R$ is a component of $w_l-x$, then the
convergence of the Fuchsian representations $\{\psi_k\circ
(h_k|_R)_*\}$ implies that any subsequence of $\{\rho_k\circ
(h_k|_R)_*\}$ has a convergent subsequence in $AH(\pi_1(R))$. This is
because, for an appropriate
choice of basepoint in $r_k(R)$, the translation distances of a
generating set of elements in $\rho_k\circ(h_k|_R)_*(\pi_1(R))$ are
bounded by the translation distances in the Fuchsian
representations (since $r_k$ are path-wise isometries).  Hence, we may
pass to a further subsequence of 
$\{\rho_k\}$, again called $\{\rho_k\}$, such that if $R$ is any
component of $w_l-x$, then $\{\rho_k\circ (h_k|_R)_*\}$ converges in
$AH(\pi_1(R))$.

Let $\{\rho_k\}$ be the subsequence and
$\{f_k:w\to w\}$ be the sequence of homeomorphisms obtained
after applying the above process $m$ times.
Each $f_k$ induces a homeomorphism $F_k:W\to W$ preserving
$\boundary M - \boundary W$, which is homotopic
to the identity on each component of $\partial W-\partial M$.
Hence, $F_k$ extends to a homotopy equivalence $\hat F_k$ of $M$ which
is equal to the identity  on the complement of a regular neighborhood of $W$.

Let $\phi_k'=(\hat F_k)_*$ and let $Z$ be the $I$-bundle over $\partial w-x$.
Thurston's relative compactness theorem and our construction of $Z\cup X$ imply
that there exists a subsequence of $\{\rho_k\}$, again denoted $\{\rho_k\}$,
such that if $M'$ is  any
component of $M-(Z\cup X)$ whose closure is not a thickened torus,
then $\{\rho_k\circ\phi_k'|_{\pi_1(M')}\}$ converges in $AH(\pi_1(M'))$.
Notice further that if $c$ is the core curve of any component of $Z$, then
$\{l_{\rho_k\circ\phi_k'}(c)\}$ converges
to a positive number $\delta(c)$. On the other hand, if $c$ is the core curve of
a component of $X$, then $\{l_{\rho_k\circ\phi_k'}(c)\}$ converges to 0.

In order to complete the proof, it is necessary to  
precompose by an additional sequence of automorphisms (and pass
to an additional subsequence)
to guarantee that the new sequence of representations converges
on every component of $M-X$ which is not a thickened torus.
We will make repeated use of the following elementary lemma,
which we state without proof.

\begin{lemma}{Dehntwist}{} Let $Q$ be a compact, hyperbolizable 3-manifold and
let $A$ be an essential annulus in $Q$ with core curve $a$.
Suppose that $\{\phi_i\}$ is
a sequence in $AH(\pi_1(Q))$ such that
if $Q'$ is any component of $Q-A$, then $\{\phi_i|_{\pi_1(Q')}\}$ converges
in $AH(\pi_1(Q'))$ and such that $\{l_{\phi_i}(a)\}$ converges to a
positive number $\delta$. Then,
there exists a subsequence $\{\phi_k\}$ of $\{\phi_i\}$ and a sequence
$\{ g_k\}$ of homeomorphisms of $Q$ (each of which is a power of
a Dehn twist about the annulus $A$) such
that $\{\phi_k\circ (g_k)_*\}$ converges in $AH(\pi_1(Q))$.
\end{lemma}

Let $Z_1,\ldots, Z_n$ be the components of $Z$.
We inductively assume that we have chosen a subsequence $\{\rho_j\}$
of $\{\rho_k\}$ and a sequence $\{\phi_j\}$ of automorphisms of
$\pi_1(M)$ such that if $M'$ is a
component of $M-(X\cup Z_l\cup\cdots\cup Z_n)$ which is not
a thickened torus, then
$\{\rho_j\circ \phi_j|_{\pi_1(M')}\}$ converges in $AH(\pi_1(M'))$.
(If $l=1$, then $\{\rho_j\}$ is the sequence $\{\rho_k\}$ obtained above
and $\{\phi_j\}=\{\phi_k'\}$.)

Let $Q_l$ be the component of $M-(X\cup Z_{l+1}\cup\cdots\cup Z_n)$
which contains $Z_{l}$. Then
$\{\rho_j\circ\phi_j|_{\pi_1(Q_l)}\}$ and $Z_l$ satisfy the hypotheses
of Lemma \ref{Dehntwist}. We may thus assume, perhaps after passing
to a further subsequence, that there is a sequence $\{ g_j\}$ of
homeomorphisms of $Q_l$ such that $\{\rho_j\circ\phi_j\circ (g_j)_*\}$
converges in $AH(\pi_1(Q_l))$. Since each $g_j$ is a power of a Dehn
twist about the annulus $Z_l$,  we may extend $\{g_j\}$ to a sequence,
again called $\{ g_j\}$,
of homeomorphisms of $M$, such that each $g_j$ is the identity map on
$M-Q_l$. Then, after  replacing $\phi_j$ by $\phi_j\circ (g_j)_*$, we
see that if $M'$ is any component of $M-(X\cup Z_{l+1}\cup\cdots\cup Z_n)$
which is not a thickened torus, then
$\{\rho_j\circ \phi_j|_{\pi_1(M')}\}$ converges in $AH(\pi_1(M'))$.

Applying the above process $n$ times gives the desired subsequence
and sequence of automorphisms. This completes the proof of
Theorem \ref{unmarked-thurston}.
\end{pf}

\section{3-manifolds which are not generalized books of I-bundles}
\label{unbook}

The main result of this section asserts
that hyperbolizable 3-manifolds with incompressible boundary
which are not generalized books of $I$-bundles have $\Lambda$-invariant
strictly less than 1.

\medskip
\par\noindent
{\bf Theorem \ref{notbook}.} {\em
If $M$ is a compact, orientable, hyperbolizable 3-manifold
with incompressible boundary
which is not a generalized book of $I$-bundles, then $\Lambda(M)<1$.}

\medskip

\begin{pf*}{Proof of Theorem \ref{notbook}}
Let $M$  be a compact, orientable, hyperbolizable 3-manifold  with incompressible boundary
which is not a generalized book of $I$-bundles.
We will assume that $\Lambda(M)=1$ and arrive at a contradiction.

If $\Lambda(M)=1$, there exists a sequence $\{ N_i\}$ in $TT(M)$ such
that $\{\lambda_0(N_i)\}$ converges to 1 where $N_i=\Hyp^3/\Gamma_i$.
Let $\rho_i:\pi_1(M)\to \PSL 2(\C)$ be a discrete faithful representation
with image $\Gamma_i$.
Let $\{\rho_j\}$, $\{\phi_j\}$, $X$,
and $x$ be as in Theorem \ref{unmarked-thurston}.

Let $M_0$ be a component of $M-X$ which contains a 
component $V$ of $M-\Sigma(M)$ whose closure
is not a solid torus or a thickened torus.
Note that the fundamental group of $V$ must be non-abelian, since
any compact hyperbolizable 3-manifold with abelian fundamental group
is a solid torus or a thickened torus.
Hence, there exists a curve in $V$ which is not homotopic
into $\Sigma(M)$ and the closure of $M_0$ cannot be a thickened torus.
Thus, Theorem \ref{unmarked-thurston} implies that
$\{\rho_j\circ\phi_j|_{\pi_1(M_0)}\}$ converges
to a discrete faithful representation $\rho:\pi_1(M_0)\to \PSL 2(\C)$
such that $\rho(\pi_1(M_0\cap X))$ is parabolic.
Let $\Gamma^0=\rho(\pi_1(M_0))$  and $N^0=\Hyp^3/\Gamma^0$.
Let $\hat\Gamma$ be a geometric limit of some subsequence of
$\{\Gamma_j^0=\rho_j(\phi_j(\pi_1(M_0)))\}$ and let $N_j^0=\Hyp^3/\Gamma_j^0$
and $\hat N=\Hyp^3/\hat\Gamma$.
Since
$1\ge \lambda_0(N_j^0)\ge\lambda_0(N_j)$
and $ \lim_{j\to\infty}\lambda_0(N_j)=1$, we see that
$\lim_{j\to\infty}  \lambda_0(N_j^0)= 1$.
Lemma \ref{lambda-semicontin} then implies that $\lambda_0(\hat N)=1$ and
hence that $\lambda_0(N^0)=1$.

Let $M_0'$ denote $M_0-{\cal N}(X)$ where ${\cal N}(X)$ is a regular
neighborhood of $X$ and let $Y$ denote the intersection of
the closure of ${\cal N}(X)$ with $M_0'$.
Since $\partial M_0'-Y$ is incompressible and the elements
of $\Gamma^0$ corresponding to $\pi_1(Y)$ are parabolic,
Lemma \ref{Bcondition} implies that $\Gamma^0$ satisfies
Bonahon's condition (B).  Therefore, $N^0$ is topologically tame.
Since, $\lambda_0(N^0)=1$, the results of \cite{canary-taylor} imply that
$\Gamma^0$ is either a Fuchsian group 
or an extended Fuchsian group. Thus $N^0$ contains
a compact submanifold $R$ which is a strong deformation retract of $N^0$ and
is an I-bundle with base surface $B$,
such that an element of $\Gamma^0$ is parabolic if and only if
it is conjugate to an element of $\pi_1(\partial B)$.
(If $Z$ is the totally geodesic hyperplane preserved by $\Gamma^0$,
then we can take $B$ to be a compact core for $Z/\Gamma^0$ and
$R$ to be a regular neighborhood of $B$.)

Since $\Gamma^0=\rho(\pi_1(M_0))$ and every element of $\Gamma^0$ corresponding
to an element of $\pi_1(Y)$ is parabolic, there is a homotopy equivalence
$h:M_0'\to R$ such that  every element of $Y$ maps into a component
of $S$ the sub-bundle over $\partial B$.
If we let $S_0$ denote the set of components of $S$ which contain
images of elements of $Y$, then $h:(M_0',Y)\to (R,S_0)$ is a homotopy
equivalence of pairs. Since every component of $\partial M_0'-Y$
is incompressible and $h$ is a homotopy equivalence of pairs,
every component of $\partial R-S_0$ is incompressible
(see Proposition 1.2 in \cite{bonahon} or section 2 of
Canary-McCullough \cite{canary-mccullough}), 
which implies that $S_0=S$. Thus, $(M_0',Y)$ is
homotopy equivalent to the $I$-pair $(R,S)$ which implies that $(M_0',Y)$ is
homeomorphic to $(R,S)$ (see corollary 5.8 in Johannson \cite{Johannson}).
This implies that $M_0'$ is an admissibly embedded essential $I$-bundle
and hence may be properly  homotoped into $\Sigma(M)$.
This however contradicts the fact that there exist curves in $V$
which are not homotopic into $\Sigma (M)$.
This contradiction completes the proof of Theorem \ref{notbook}.
\end{pf*}

\section{Remarks and Conjectures}
\label{remarks}

{\bf 1.} The most natural analogue of the Gromov norm is the invariant
$V(M)$ which is defined to be the infimum of the volumes of the convex
cores of hyperbolic manifolds homeomorphic to the interior of $M$.
One would make the following conjecture in the spirit of the paper.

\medskip

{\bf Conjecture:} {\em $V(M)=0$ if and only if every component of
$M$'s incompressible core is a generalized book of $I$-bundles.}

\medskip

It seems likely that  an argument similar to that in section \ref{books}
would give that $V(M)=0$ whenever $M$ is a generalized book of $I$-bundles.
Whereas an argument similar to that in section \ref{unbook}, along
with work of Taylor \cite{taylor}, should give that $V(M)>0$ if
$M$ has incompressible boundary but is not a generalized book of $I$-bundles.
If $M$ is obtained from generalized books of $I$-bundles by adding
1-handles, then one needs to show that $V(M)=0$.
However, as the direct analogues for volumes of
the results of Patterson for $\lambda_0$, used 
in Section \ref{compress}, are false,
one would need to find a more explicit proof for this case.

\bigskip

{\bf 2.}
The work of Canary \cite{CanaryJDG} and Burger-Canary \cite{burger-canary}
exhibits relationships between $V(M)$ and $\Lambda(M)$.
The work of \cite{CanaryJDG} gives the following upper bound for all $M$,
$$\Lambda(M)\le { 4\pi|\chi(\partial M)|\over V(M)}$$
(where $\chi(\partial M)$ denotes the Euler characteristic of $\partial M$),
while the work of \cite{burger-canary} shows that there exist constants
$A>0$ and $B>0$ such that if $M$ has incompressible boundary
then
$$\Lambda(M)\ge {A\over (V(M)+B|\chi(\partial M)|)^2}.$$

\bigskip

{\bf 3.} If $M$ is acylindrical then there is a unique
hyperbolic manifold $N$ whose convex core has totally geodesic
boundary and is homeomorphic to $M$.
One expects that $D(M)=d(N)$, $\Lambda(M)=\lambda_0(N)$, and
that $V(M)$ is the volume $vol(C(N))$ of the convex core $C(N)$ of $N$.
In fact, Bonahon \cite{bonahon2} has shown that $N$ is a local minimum for the
volume (of the convex core) function on the space $GF(M)$
of geometrically finite hyperbolic 3-manifolds homeomorphic to
the interior of $M$.  More generally, if $M$ has incompressible
boundary and $M_i$ is a component of $M-\Sigma(M)$ which is
not homeomorphic to a solid torus, then there exists
a unique hyperbolic 3-manifold $N_i$ such that the convex core
$C(N_i)$ has  totally geodesic boundary
and $C(N_i)$ is homeomorphic to $M_i$.

\medskip
\par\noindent
{\bf Conjecture:} {\em If $M$ has incompressible boundary,
$\{ M_1,\ldots,M_n\}$ are the components of $M-\Sigma(M)$
which aren't homeomorphic to solid tori, and $N_i$ are as above,
then $V(M)=\sum_{i=1}^n vol(C(N_i))$ and $\Lambda(M)=min \{ \lambda_0(N_i)\}$.}

\medskip

A positive solution to the conjecture above would imply that $\Lambda$
is a homotopy invariant, since 
the incompressible cores of
homotopy-equivalent hyperbolic 3-manifolds are homotopy-equivalent, 
and Johannson's theorem \cite{Johannson} implies
that the complements of the characteristic submanifolds of
homotopy-equivalent hyperbolizable 3-manifolds with incompressible
boundary are homeomorphic. 

We note that $\sum_{i=1}^n vol(C(N_i))$ is equal to half the Gromov
norm of the double of $M$.
One can generalize this conjecture
to obtain a similar conjecture for arbitrary compact hyperbolizable  3-manifolds.

\bigskip

{\bf 4.} It would be interesting to know more about the distribution 
of the set of values assumed by the invariant $\Lambda$.
In this remark we show that the set of values accumulates at 0.

Let $M_i$ be a compact hyperbolizable 3-manifold whose
boundary is incompressible and has $i$ toroidal
boundary components and two genus two boundary components.
One may obtain such manifolds by removing suitably chosen collections
of boundary-parallel curves from 
a product $S\times[0,1]$ where $S$ is a surface of genus 2.
We will show that $V(M_i)\to \infty$, by observing that all but a
bounded number of the cusps contribute a definite amount to the volume
of the convex core of {\em any} hyperbolic structure on $M_i$.   

If $x\in N$, then $inj_N(x)$ denotes the injectivity radius of the
$N$ at the point $x$.
The Margulis lemma asserts that there exists a constant ${\cal M}_3$
such that if $\epsilon<{\cal M}_3$ and $N$ is a complete hyperbolic
3-manifold, then every non-compact component
of $N_{thin(\epsilon)}=\{ x\in N|inj_N(x)<\epsilon\}$ is the quotient
of a horoball by a group of parabolic transformations fixing the horoball.
Moreover, for all ${\cal M}_3>\epsilon >0$, there exists $D(\epsilon)>0$
such that any non-compact component of $N_{thin(\epsilon)}$ has 
volume at least $D(\epsilon)$.
There also exists $K(\epsilon)>0$ such that
at most $K(\epsilon)|\chi(\partial C(N_i))|$ components
of $N_{thin(\epsilon)}$ intersect $\partial C(N)$.
If $N_i\in TT(M_i)$  then there are at least $i$ non-compact
components of $(N_i)_{thin(\epsilon)}$ (one for each toroidal component
of $\partial M_i$.)
Since $|\chi(\partial C(N_i))|\le 4$ and every component of
$(N_i)_{thin(\epsilon)}$ intersects $C(N_i)$,
it follows that 
$$vol(C(N_i)) >(i-4K(\epsilon)) D(\epsilon).$$
Hence, $V(M_i)\ge (i-4K(\epsilon))D(\epsilon)$, so $V(M_i)$ converges to
infinity and $\Lambda(M_i)\le {16\pi\over V(M_i)}$ converges to 0.
On the other hand, $\Lambda(M_i)\ne 0$, by Proposition \ref{whenDis2}.
Thus $0$ is an accumulation point of the set of $\Lambda$ values.

\bigskip
{\bf 5.} If $\Gamma$ is a quasi-Fuchsian group, denote by $K(\Gamma)$
the minimal $K$ for which there exists a $K$-quasiconformal
map conjugating $\Gamma$ to a Fuchsian group.
One may use the methods of section \ref{books} to construct a
sequence $\{\Gamma_j\}$ of quasiconformally conjugate quasi-Fuchsian groups
such that the Hausdorff dimensions of the limit sets of the
$\Gamma_j$ converge  to 1,
but $\{K(\Gamma_j)\}$ converges to $\infty$. We will discuss this more
fully in a future note.

\end{document}